\documentclass{amsart}
\usepackage{amsthm}
\usepackage{amsmath}

\title{Degenerations of Representations and Thin Triangles.}  
\author{Daryl Cooper}
\thanks{Research supported in part by the NSF circa 1992.}

\newtheorem{theorem}{Theorem}[section]
\newtheorem{lemma}[theorem]{Lemma}
\newtheorem{corollary}[theorem]{Corollary}

\newtheorem{proposition}[theorem]{Proposition}
\newtheorem{addendum}[theorem]{Addendum}

\newcommand\Hom{\rm{Hom}}
\newcommand\Isom{\rm{Isom}}
\newcommand\C{\Gamma}
\newcommand\s{\sigma}

\newcommand\od{\overline{d}}

\newcommand\e{\epsilon}

\newcommand\io{\iota}

\newcommand\nb{\newline$\bullet\ \ \ $}
\newcommand\PLF{{\mathbb P{\rm LF}}}
\newcommand\LF{{{\rm LF}}}
\newcommand\demo{\begin{proof}}
\newcommand\SL{\operatorname{SL}}

\newcommand\RR{\mathbb R}
\newcommand\PP{\mathbb P}
\newcommand\Hcal{\mathcal H}
\newcommand\CC{\mathbb C}
\newcommand\ZZ{\mathbb Z}
\newcommand\HH{\mathbb H}
\newcommand\Newt{\operatorname{Newt}}

\title{Degenerations of Representations and Thin Triangles.}  

\author{Daryl Cooper}
\address{Math Dept,
University of California ,
Santa Barbara,
CA. 93106, 
U.S.A.}
\curraddr{}
\email{}
\thanks{Research supported in part by the NSF circa 1992.}


\date{}

\begin{document}
\begin{abstract} There is a compactification of the space of representations of a
finitely generated group into the groups of isometries of all spaces with $\Delta$-thin
triangles. The ideal points are actions on ${\mathbb R}$-trees. It is a geometric reformulation
and extension of the Culler-Morgan-Shalen theory concerning limits of 
representations into $\SL_2{\CC}$ and more generally ${\rm SO}(n,1)$.
\end{abstract}

\maketitle

This paper was written and circulated in the early 90's,  but never published.
Much has happened since then, but we have not attempted to update it. 

Let $G$ be a
finitely generated group. A map $L:G\to\RR$ is called a {\em length function}.
The space of  length functions is $\LF(G)=\RR^G$ with
the product topology. This is compactified by the
space, $\PLF(G)=\PP(\RR^G)$,  of projectivized length functions.

Suppose $\Hcal$ is a geodesic space with $\Delta$-thin triangles and isometry group $\Isom(\Hcal)$. 
 A homomorphism $\rho:G\to\Isom(\Hcal)$ determines a
length function $L(\rho):G\longrightarrow {\mathbb R}$  
$$L(\rho)(g)=\inf\{d_\Hcal(x,\rho(g)x)\ :\ x\in\Hcal\}$$
The collection of all such $\rho$ is denoted $\Hom_{\Delta}(G)$, and there is a map
$$L:\Hom_{\Delta}(G)\longrightarrow \LF(G)$$ The  image of $L$ has compact closure, and limit points in $\PLF(G)$ are projective classes of
length functions of certain actions of $G$ on ${\mathbb R}$-trees. This is because a limit of rescaled
$\Delta$-thin spaces  is an $\mathbb R$-tree.

This may be viewed as a generalization of the compactification of spaces of characters of
representations into the Lie groups ${\rm SO}(n,1)$ used by Culler, Morgan and Shalen
\cite{CS1,CS2,CGLS,MS} by actions on $\mathbb R$\--trees. 
This idea was used by
Bestvina \cite{Be} to provide new proofs of some results of Morgan and Shalen \cite{MS}. In
particular he applied a geometric method to obtain this compactification for the subspace of
discrete faithful representations into ${\rm O}(n,1)$. This paper removes the hypotheses
of discrete and faithful, and applies the technique to a wider class of groups than ${\rm O}(n,1)$.

Next we give sufficient conditions on a length function to ensure it is the length
function of an action of $G$ on a {\em simplicial} tree. This is used to recover some results of
Culler and Shalen \cite{CS1} on compactifications of curves of characters by actions on
simplicial trees. The standard approach to this uses the the Tits\--Bass\--Serre
approach \cite{Se} to trees via ${\rm SL}_2\mathcal F,$ where $\mathcal F$ is a field with a discrete rank-1
valuation. We have reduced the use of valuations to showing that a certain ${\mathbb R}$-tree
has a minimal invariant subtree that is simplicial. 
 
In what follows $\Hcal$ is   a geodesic space with thin triangles.
The starting point in section 1 is the observation that if
 the vertices of a finite connected
graph in $\mathcal H$  are moved, so that the diameter of the
graph  goes to infinity, and if the induced metric on the graph is rescaled to keep the diameter
bounded, then the metric on the graph converges to a pseudo-metric for which the associated
metric space is a tree. 

Let $K$ be the complete graph with vertex set $G$. 
Given a homomorphism $\rho:G\longrightarrow \Isom(\Hcal)$ there is a $G$-equivariant map of $f:K\to \Hcal$.
Each edge maps linearly.
 The map is determined by the choice of a single point $x\in\Hcal$ that  is a
modified version of the {\it center} of $\rho$ introduced by Bestvina \cite{Be}. In some sense $\rho$ {\em operates}
near $x$.

Given a sequence $\rho_n:G\to\Isom(\Hcal_n)$ of homomorphisms for which the translation length of some element goes to
infinity,  there is a subsequence such that the rescaled pullback metrics on $K$ converge to an $\mathbb R$-tree, $\Gamma,$ on
which $G$ acts by isometries. This leads to the compactification result. We then show that, if
the translation length of every element of $G$ acting on $\Gamma$ is an integer, and if the
action is irreducible, then $\Gamma$ has an invariant {\em simplicial} subtree.  This is done in section 2.

Section 3 discusses how spaces of representations can be compactified by mapping the
representation space into the space of length functions on the group and compactifying this
latter space projectively. This is how Thurston compactified Teichmuller space.
The argument up
to this point is geometric and combinatorial. 

In section 4, we specialize to the action of $\SL_2\CC$ by isometries on $\HH^3$ .
Valuations make their appearance to show that, if the $\rho_n$  all
lie on an {\it algebraic curve}  $C\subset \Hom(G,\SL_2{\CC})$ then, after rescaling the limit $\mathbb R$-tree,
 the translation length of every
element is an integer. 
The idea is that   the translation length of an element $g\in G$ is approximately $\log (1+|p_g|)$ where
 $p_g=\operatorname{trace}(\rho g)$ is a function of $\rho\in C$.

If $p(z)$ is a polynomial of degree $d$  in {\em one variable $z$}, 
then $\log (1+|p(z)|)\sim d \log |z|$ for  $z$ large.  For each end $\epsilon$ of $C$
there is a valuation $\nu_{\epsilon}$ on the function field of $C$,
 that generalizes the degree of a polynomial. It can be scaled to take integer values. 
If $f:C\to\CC$ is  an algebraic function and $z\in C$ moves out into  $\epsilon$,   
then $f(z)$ behaves  like a rational function of one variable, with a degree of growth given by $\nu_{\epsilon}(f)$.
This culminates in the main theorem (\ref{mainthm}).

These ideas were  surely known to Culler, Morgan and Shalen, and
motivated the search for a more algebraic formulation. However we believe the
geometric point of view deserves to be more widely known. 

\section{Rescaling Graphs in Negatively Curved Spaces}

The results of this section apply to geodesic spaces satisfying the {\it thin triangles}
property introduced by Gromov in \cite{Gr}. The main result of this section is (\ref{Rtree}) which
states that the limit of certain rescaled metric spaces is an ${\mathbb R}$-tree. 

Suppose $(X,d)$ is  a  metric space. The
{\em length} of a path $\gamma:[a,b]\longrightarrow X$ is defined to be the supremum over finite
dissections of $[a,b]$ of the sum of the distances between endpoints of the dissecting
subintervals. A {\em path-metric} is a metric with the property that the distance between two
points is the infimum of the lengths of paths connecting the points. Thus a Riemannian manifold is a
path-metric space. 

A {\em metric simplicial tree} is a path-metric space that is also a simplicial  tree.
 An {\em ${\mathbb R}$-tree} is a path-metric space, $\Gamma$, with the
property that given any two points $a,b\in\Gamma$ there is a unique arc $\gamma$ in $\Gamma$ with endpoints $a,b$
and $\gamma$ is isometric to an interval in $\RR$, see \cite{CM} or \cite{Sh}. Thus every metric
simplicial tree is an ${\mathbb R}$-tree.

Suppose $X$ is a path metric space. A {\em geodesic}  is an arc in $X$ with the property that the
length of the arc equals the distance in $X$ between its endpoints. In a
simply-connected complete Riemannian manifold of non-positive curvature this is equivalent to
the usual definition of geodesic. A path-metric space in which every pair of points is contained in a geodesic is called a {\em geodesic space.} 
Given two points $a,b\in X$ then $[a,b]$ denotes some choice of geodesic
with endpoints $a,b$. A {\em degenerate geodesic} $[a,a]$ consists of a single
point. 

A {\em triangle} in a geodesic space consists of three geodesics $[a,b],[b,c],[c,a]$
called the {\em sides} or {\em edges} of the triangle. Given $\Delta>0$ a triangle is called $\Delta$-{\em thin} if
each side is contained in the $\Delta$-neighborhood of the union of the other two sides. A geodesic space is called {\em negatively
curved} or {\em $\Delta$-hyperbolic} of {\em Gromov hyperbolic} if there is $\Delta>0$ such that every triangle is
$\Delta$-thin. 

It is clear that simplicial trees, and more generally ${\mathbb R}$-trees, have
$\Delta$-thin triangles for all $\Delta>0$. The term {\em hyperbolic space}  refers
to ${\mathbb H}^n$, which is the complete simply connected Riemannian $n$-manifold with constant
sectional curvature $-1$.
It is well known that triangles in ${\mathbb H}^n$ are $\Delta$-thin with
$\Delta=\log\left(1+\sqrt{2}\right)$.

\begin{proposition}\label{thinpoly} Suppose $(\Hcal,d_\Hcal)$ is a geodesic space with $\Delta$-thin
triangles, and $P$ is a polygon in $\Hcal$ that
has $n$ geodesic sides. Then every side of $P$ lies within a distance $(n-2)\Delta$ of the union
of the other sides of $P$\end{proposition} 
\demo The polygon $P$ may be triangulated using $n-2$ triangles by coning from a vertex of $P,$ and
the result follows from the thin triangles property. \end{proof}

From now on we will assume that $(\Hcal,d_\Hcal)$ is a geodesic space with $\Delta$-thin
triangles. Often we will write $d(x,y)$, instead of $d_\Hcal(x,y)$.
Let $K$ be a {\em connected graph}, in other words a connected 1-dimensional simplicial complex. A
map $f:K\longrightarrow\Hcal$ is called a {\it straight map} if, for each edge $e$ of $K$, the restriction of $f$ to $e$ is {\em linear} in the following sense.
Since $e$ is a $1$-simplex, we may identify $e$ with $[0,1]$ by an affine map, and {\em linear} means
$d_{\Hcal}(f(0),f(t))=t\cdot d_{\Hcal}(f(0),f(1))$.

A {\em pseudo-metric} satisfies the same conditions as a metric, except that the
distance between two distinct points might be zero. There is a pseudo-metric, $f^*d_{\Hcal}$, on $K$
defined by pulling back the metric on $\Hcal$ that is given by
$$f^*d_{\Hcal}(x,y)=d_\Hcal(fx,fy)$$ A
pseudo-metric space $(K,d)$ determines an equivalence relation on $K$ given by $x\sim y$ if
$d(x,y)=0$. Then $d$ induces a metric, also denoted $d,$ on $K/\sim$. We call $(K/\sim,d)$ the
metric space {\it associated } to $(K,d)$. Clearly $f$ induces an isometry of $(K/\sim,d)$
with the subset $f(K)$ of $\Hcal$. 
However the associated metric for $f^*d_{\Hcal}$ is usually not a
path-metric because there may be no geodesic in  $f(K)$ that connects some pair of points in $f(K)$. 

\begin{proposition}\label{simplicialtree} Given $\Delta>0$, and a finite connected graph $K$,  suppose that $f_n:K\longrightarrow
\Hcal_n$ is a sequence of straight maps  into  geodesic
spaces $\Hcal_n$ with $\Delta$-thin triangles. 
Suppose that $\lambda_n\to0$, and  $\lambda_n\cdot f_n^*d_{\Hcal_n}$ converges pointwise on $K$ to a pseudo-metric
$\od$. 
Then the metric space associated to $(K,\od)$ is a metric simplicial tree\end{proposition}
\demo Let
$(\C,\od)$ be the metric space associated to $\od$.  Since $f_n$ is linear on each edge it is clear that subset
of $\Gamma$ corresponding  to an edge of $K$ is an interval, possibly with zero length.  Hence $\gamma$ is the union of finitely many
metric intervals.

The proof proceeds by induction on the number of edges of $K$. It suffices to prove the result for $K^+ =K\cup\beta$
where $\beta$ is an edge, and $K\cap\beta$ is one or both endpoints of $\beta$. Let $\beta_-\subset\beta$
 consist of all points $x\in\beta$ with the property that there is  $y\in K$ and $\od(x,y)=0$.
 We claim that $\beta_-$ is a closed subinterval of $\beta$. For notational simplicity we assume that all the spaces $\Hcal_n$ are the same.

To prove the claim, suppose $x',y'\in\beta_-$. There are $x,y\in K$ with $\od(x,x')=0$ with $\od(y,y')=0$. Since $K$ is connected, we may choose an arc $\ell$ in
$K$ that contains $x$ and $y$. Then $\ell$ is the union of some finite number, $m,$ of intervals each contained in
some edge of $K$. Thus $f_n(\ell)$ is a polygonal path in $\Hcal$ with endpoints $f_nx$ and  $f_ny$ that is the union of
$m$ geodesic segments, where $m$ does not depend on $n$.  

Let $\gamma_n$ be the geodesic segment in $\Hcal$ with
endpoints $f_nx$ and $f_ny$. Every point of $\gamma_n$ lies within a distance $D=(m-1)\Delta$ of $f_n(\ell)$. This
follows from  (\ref{thinpoly}) applied to the polygon $\gamma_n\cup f_n(\ell)$. Now
$f_n\beta_-$ is a geodesic arc in $\Hcal$ with endpoints $f_nx'$ and $f_ny'$. Since $\od(x,x')=0$
it follows that $\lambda_nd_\Hcal(f_nx',f_nx)\to 0$ as $n\to\infty$. Similarly
$\lambda_nd_\Hcal(f_ny',f_ny)\to 0$. 

Suppose $[a_1,a_2]$ and $[b_1,b_2]$  are geodesic arcs in $\Hcal$ and $d_{\Hcal}(a_i,b_i)\le M$ for $i=1,2$.
By applying (\ref{thinpoly}) to the quadrilateral $[a_1,a_2]\cup[a_2,b_2]\cup[b_1,b_2]\cup[b_2,a_1]$ it follows
 $[a_1,a_2]$ is contained in a neighborhood of size $M+2\Delta$ of  $[b_1,b_2]$.
 
We
apply this to the arcs $\gamma_n$ and $f_n(\beta_-)$. It follows that given $z'$ on $\beta_-$
there is $z_n$ on $\gamma_n$ with $d_\Hcal(f_nz',z_n)\le M+2\Delta$. Since the point $z_n$
on $\gamma_n$ is within a $d_\Hcal$-distance $D$ of some point  on $f_n(\ell)$ it
follows that $d_n(z',\ell)\le M+2\Delta+D$. Thus $\od(z',\ell)=\lim \lambda_n
d_n(z',\ell)=0$.  This proves that $\beta_-$ is an interval. It is easy to show that
$\beta_-$ is closed, establishing the claim.

Since $\beta$ is attached to $K$ at one endpoint, if $\beta-\beta_-$ is not empty, then it is an
interval. If $\beta_-=\beta$ then it is clear that the tree formed by adding $\beta$ to $K$ is
the same as the tree formed from just $K$. Let $\gamma=[a',b']$ be the closure of
$\beta-\beta_-$. One of the endpoints, say $a',$ of $\gamma$ is a $\od$-distance of zero from
some point, say $a,$ of $K$. Thus $\od(a,a')=0$. The topological space
$\Gamma\equiv(K\cup\gamma)/\sim$ is the identification space formed from the union $K/\sim$ with
$\gamma$ by identifying $a'$ with $a$. Therefore the result is homeomorphic to a simplicial tree.
It remains to show that the metric induced by $\od$ makes $\Gamma$ into a metric tree.

Let $x\in \gamma$ and $y\in K$ be arbitrary points then there are arcs $[y,a]\subset K/\sim$
and $[a',x]\subset \gamma$. The union of these arcs, $[y,a]\cup[a',x],$ is the unique arc in
$\gamma$ from $x$ to $y$. The proof is completed by the claim that 
$\od(y,x)=\od(y,a)+\od(a,x)$.

Suppose that, given $\mu>0$ that for all sufficiently large $n$ that there is a point, $z,$ on
$[f_nx,f_ny]$ with $d_\Hcal(f_na,z)\le \mu\lambda_n^{-1}+2\Delta$. Then 
$$\begin{array}{rcl}
 &    &|\ d_n(x,y)-\left( d_n(x,a)+d_n(a,y)\right)\ | \\
& \le & |\ d_\Hcal(f_nx,z)+d_\Hcal(z,f_ny)-d_\Hcal(f_nx,f_na)-d_\Hcal(f_na,f_ny)\ |\\  
& \le & |\ d_\Hcal(f_nx,z)-d_\Hcal(f_nx,f_na) +d_\Hcal(z,f_ny)-d_\Hcal(f_na,f_ny)\ | \\ 
& \le & |\ d_\Hcal(z,f_na) +d_\Hcal(z,f_na)\ | \\ 
& \le & 2(\mu\lambda_n^{-1}+2\Delta)\end{array}$$
Multipying both sides by $\lambda_n$ and using that $\lambda_n\to0,$ the claim follows in this case.

We consider the thin triangle with sides $[f_nx,f_na]$, $[f_nx,f_na']$ and $[f_na,f_na']$. Since
$\od(a,a')=0$ it follows that given $\epsilon>0$ for sufficiently large $n$ we have
$d_n(a,a')\le\epsilon\lambda_n^{-1}$. Hence for all points $q_n$ on $[f_nx,f_na]$
$$d_\Hcal(q_n,[f_nx,f_na'])\le\epsilon\lambda_n^{-1}+\Delta$$

We may assume that there is some $\mu>0$ such that for arbitrarily large $n$
that $d_\Hcal(f_na,[f_nx,f_ny])> \mu\lambda_n^{-1}+2\Delta$. Define $p_n$ to be the point on
$[y,a]$ with $d_n(p_n,a)=\mu\lambda_n^{-1}$ and $q_n$ to be the point on
$[f_nx,f_na]$ with $d_\Hcal(q_n,f_na)=\mu\lambda_n^{-1}$. Now:
$$\begin{array}{rcl}
d_\Hcal(f_np_n,[f_ny,f_nx]) &\ge& d_\Hcal(f_na,[f_ny,f_nx])-d_\Hcal(f_na,f_np_n)\\
&\ge& \mu\lambda_n^{-1}+2\Delta-\mu\lambda_n^{-1}\\
&=&2\Delta\end{array}$$
Then thin triangles applied to the triangle with vertices $f_na,f_nx,f_y$ gives $$d_{\mathcal
H}(f_np_n,[f_na,f_nx])\le\Delta$$ Let $r_n$ be a point on $[f_na,f_nx]$ with $d_{\mathcal
H}(f_np_n,r_n)\le\Delta$. Then $$|\ d_\Hcal(f_na,r_n)-d_\Hcal(f_na,f_np_n)\ | \le \Delta$$
Since $r_n$ and $q_n$ are both on $[f_na,f_nx]$ and $d_n(a,p_n)=d_\Hcal(f_na,q_n)$ it follows that
$d_\Hcal(r_n,q_n)\le\Delta$. Thus $$d_\Hcal(f_np_n,q_n) \le d_\Hcal(f_np_n,r_n)+d_{\mathcal
H}(r_n,q_n)\le2\Delta$$ Combining this with the bound on $d_\Hcal(q_n,[f_nx,f_na'])$ gives
$$d_\Hcal(f_np_n,[f_nx,f_na'])\le d_\Hcal(f_np_n,q_n)+d_\Hcal(q_n,[f_nx,f_na']) \le
2\Delta+\epsilon\lambda_n^{-1}+\Delta$$ As $n\to\infty$ there is a subsequence of $p_n$ which
converges to a point p in $[a,y]$. We get $\od(p,[x,a'])\le\epsilon$. Since $\epsilon>0$ is
arbitrary this gives $\od(p,[x,a'])=0$. Now $\od(p,a)=\mu>0$ so $p\ne a$ but this contradicts that
$a$ is the only point on $[a,y]$ zero distance from $[x,a']$. This proves the claim\end{proof}

\begin{lemma} Suppose that $(\Gamma,d)$ is a metric space which is the union of an increasing
sequence of subspaces $T_n,$ each of which is a metric simplicial tree. Then $\Gamma$ is an
${\mathbb R}$-tree\end{lemma}
\demo Suppose that $x,y$ are two distinct points in $\Gamma$ then there is some $T_m$ which
contains both of them. There is a unique arc $[x,y]$ in $T_m$. Suppose that there is a
different arc $\alpha$ in $\Gamma$ with the same endpoints. The intersection of $\alpha$ with
$[x,y]$ is closed. We may assume that $[x,y]$ intersects $\alpha$ only in the points $x,y$.
For otherwise we may replace $[x,y]$ by a suitable subarc.

Define a nearest-point projection $$\pi:\Gamma\longrightarrow[x,y]$$ as follows. 
Given a point $z$ in $\Gamma$ there is an $n>m$ such that $T_n$ contains $z$. Define $\pi z$ to
be the unique point on $[x,y]$ which is closest to $z$ in $T_n$. Since the trees $T_n$ are an
increasing sequence, it is clear that this definition of $\pi z$ is independent of $n$. We
claim that $\pi$  is continuous. For if $d_{\Gamma}(z,z')<\epsilon$ then if $n$ is chosen large
enough that $T_n$ contains both $z$ and $z'$ then $d_{T_n}(z,z')<\epsilon$. But the
restriction of $\pi$ to the simplicial tree $T_n$ clearly does not increase distance thus
$d_{T_n}(\pi z,\pi z')<\epsilon$. This proves the claim.

We claim that $\pi\alpha=x$. Let $\beta$ be the maximal subarc of $\alpha$ containing $x$ such
that $\pi \beta=x$.  Since $\pi x=x$ it follows that $\beta$ is not empty. Since $\pi$ is
continuous it follows that $\beta$ is closed. Writing $\beta=[x,z]$ choose $n>m$ with $T_n$ 
containing $z$. Assuming $z\ne x,y$ then $d_{T_n}(z,\pi z)=d_{T_n}(z,[x,y])=\delta>0$. Now $\pi y=y$
so $\beta\ne\alpha$ thus we may choose $z'$ a point in $\alpha-\beta$ with $\pi z'\ne x$. Now choose
$n>m$ such that $T_n$ contains both $z,z'$. By choice of $z'$ we have that $[z,z']=[z,\pi z]\cup[\pi
z,\pi z']\cup[\pi z',z']$. Therefore $d_{T_n}(z,z')\ge d_{T_n}(z,\pi z)=\delta$. Since
$z'$ can be arbitrarily close to $z$ this is a contradiction. Hence $z=x$.

It follows that there is a unique arc in $\Gamma$ between any pair of distinct points. Since
there is some $T_n$ containing both these points, this arc is metrized as an interval. Hence
$\Gamma$ is an ${\mathbb R}$-tree\end{proof}

\noindent Combining these two results together gives:

\begin{theorem}\label{Rtree}Suppose that $\Delta>0$, and $\Hcal_n$
is a sequence of geodesic spaces with $\Delta$-thin triangles, and $K$ is a connected graph with countable many edges. 
Suppose that $f_n:K\to\mathcal H_n$ is a
sequence of straight maps. 

Suppose that $\lambda_n\to0$, and  $\lambda_n\cdot f_n^*d_{\Hcal_n}$ converges pointwise on $K$ to a pseudo-metric
$\od$. 
Then the metric space, $\Gamma,$ associated to $(K,\od)$ is an ${\mathbb
R}$-tree\end{theorem}

\section{Limiting Actions on ${\mathbb R}$-trees.}
In this section $G$ is a finitely generated group, and we  study actions of  $G$  by isometries on a
 geodesic space $\Hcal$ with $\Delta$-thin triangles. A sequence of homomorphism ({\em representations}) $\rho_n:G\longrightarrow \Isom(\Hcal)$
 {\em blows up} if there is $\alpha\in G$ such that
$\lim_{n\to\infty}t(\rho_n\alpha)=\infty$, see \cite{CL}. If this case, after rescaling, 
there is  a limiting action of $G$ on an ${\mathbb R}$-tree $\Gamma$. With extra hypotheses, $\Gamma$ contains a  $G$-invariant 
simplicial subtree.

Suppose that $(X,d)$ is a metric space and $\tau$ is an isometry of this space then 
the {\em translation length} of $\tau$ is
$$t(\tau)=\inf\{d(x,\tau x):\ x\in X\ \}$$
This defines a map $t:\Isom(X)\longrightarrow\mathbb R$. 

Let $K\equiv K(G)$ be the complete graph with $G$ 
as vertex set. This is a graph with one edge between every pair of vertices. Choose a finite
generating set $S$ for $G$. Given $x$ in $\Hcal$, and a  representation
$\rho:G\longrightarrow \Isom(\Hcal)$, define  $$r_S(x)=\max\{d(x,(\rho\alpha)x):\ \alpha\in S\ \}$$ An
{\it approximate center} of  $\rho$ is a choice of $x$ in $\Hcal$ for
which $r_S(x)\le r_S(y)+1$ for all $y\in\mathcal H$, compare
\cite{Be}. 

For example, if
$G$ is the fundamental group of a compact convex hyperbolic surface with boundary, the
developing map identifies the basepoint of the universal cover with some point in the hyperbolic
plane. This is an approximate center. Intuitively   the representation {\em operates near here}.
 It is used to overcome the problems associated with changing a
representation by conjugacy.
 
Choose an approximate center $x_n$
for $\rho_n$ and define $\io_n:K\longrightarrow\Hcal$ on the vertices of $K=K(G)$ by:
$$\io_n(\alpha)=(\rho_n\alpha)x_n$$ Then extend $\io_n$ linearly over each edge of $K$ so that $\io_n$ is
a straight map. Let $d_n=\iota_n^*d_{\Hcal}$ be the pullback pseudo\--metric on $K$. Now $G$ acts on the vertices of
$K$ by multiplication on the left. Extend this action linearly over the edges of $K$ to get an
action $\s:G\longrightarrow \Isom(K,d_n)$ by isometries on the pseudo\--metric space $(K,d_n)$
It is clear that $\io_n$ is a $G$\--equivariant
isometry from $K$ onto its image.

Define the {\it rescaling factor } $$\lambda_n=\left(\max\{1,d_n(1_G,\beta):\ \beta\in S\
\}\right)^{-1}$$ and the {\it rescaled pseudo-metric } by $\tilde d_n=\lambda_nd_n$. Since $\rho_n$ is
blowing up it is easy to see that $\lambda_n\to0$.  

\begin{proposition}\label{blowupgivesRtree} Suppose that $G$ is a finitely generated group and
$\rho_n:G\to \Isom(\Hcal)$ is a sequence of representations which blows up. Let
$\tilde{d}_n$ be as above. Then there is a subsequence 
of $\tilde{d}_n$ that converges pointwise to a pseudo-metric $\od$ on $K$. Furthermore, the
metric space associated to $(K,\od)$ is an ${\mathbb R}$-tree\end{proposition}
\demo Define $K_S$ to be the subgraph of $K$ which is the Cayley graph
corresponding to the generating set $S$. Thus for each $s\in S$ there is an edge in $K_S$
connecting $g$ to $gs$. This condition is preserved by multiplication on the left, hence
$K_S$ is preserved by $G$. Thus $K_S$ is a connected graph on which $G$ acts. 

The rescaling
factor was chosen so that the $\tilde d_{n}$-length of every edge in $K_S$ adjacent to the
identity is at most $1$. Since $G$ acts transitively on vertices, and by $\tilde
d_{n}$-isometries, it follows that the $\tilde d_{n}$-length of every edge in $K_S$ is at
most $1$. Hence, for every pair of vertices $x,y$ in $K_S$ we have that $\tilde
d_{n}(x,y)$ is bounded by the number of edges in a path from $x$ to $y$. Since $\tilde d_{n}$
metrizes every edge of $K$ as an interval, it follows that for every pair of points $x,y$ in
$K$ there is a bound on $\tilde d_{n}(x,y)$

Choose a countable dense subset, $D,$ of $K$. For $x,y$ in $D$ we have that $\tilde
d_{n}(x,y)$ is bounded independently of $n$. Hence there is a subsequence, $n_i,$ so that for
all $x,y$ in $D$ then $\tilde d_{n_i}(x,y)$ converges. We claim for all $x,x'$ in $K$ that
$\tilde d_{n_i}(x,x')$ converges. The pointwise-limit is clearly a pseudo-metric, $\tilde{d},$
on $K$. In what follows we assume the sequence has been replaced by the subsequence.

To prove the claim, consider an edge $e$ of $K,$ then the map $\io_n$ is linear on $e$. The
limit pseudo-metric restricted to $e$ metrizes $e$ as a closed interval. Given a point, $x,$
on $e$ we may choose $x_1,x_2$ in $D\cap e$ with $x$ between them and with the limiting
distance between $x_1$ and $x_2$ very small so that for all $n$ sufficiently large
$\tilde{d}_n(x,x_1)<\epsilon$. Given another point $x'$ on some other edge $e'$ we may
similarly  choose $x_1',x_2'$ with limiting distance between them $\epsilon'$. Then for $n$
sufficiently large:
 $$\begin{array}{rcl}
|\tilde{d}_n(x,x')-\tilde{d}_n(x_1,x'_1)| &\le& \tilde{d}_n(x,x_1)+\tilde{d}_n(x'_1,x')\\
&\le& \tilde{d}_n(x_2,x_1)+\tilde{d}_n(x'_1,x'_2)\\
&\le&2\epsilon\end{array}$$
Since $\tilde d_n(x_1,x_1')$ converges it follows that
$\tilde d_{n}(x,x')$ converges as $n\to\infty$. It follows from (\ref{Rtree}) that the metric
space associated to $(K,\od)$ is an ${\mathbb R}$-tree\end{proof}

 We now describe isometries of $\RR$ trees, see  \cite{CM}
for a more complete treatment.
 The {\em core} of isometry, $\tau$, of an ${\mathbb R}$-tree $\Gamma$ 
is the set of points in $\Gamma$ which are moved the minimal distance by $\tau:$
 $$C(\tau)=\{\ x\in \Gamma\ :\ d_{\Gamma}(x,\tau x)=t(\tau)\ \}$$  
It is clear that this a closed convex set (ie $x,y\in {\rm core}(\tau)$ implies $[x,y]$ is in
${\rm core}(\tau)$) which is invariant under $\tau$. An isometry, $\tau,$ of $\Gamma$ is called {\em
elliptic} if $\tau$ has a fixed point. It is called {\em hyperbolic} if $t(\tau)>0$ and the core
of $\tau$ is isometric to a line. This line is called the {\em axis} of $\tau$. 

\begin{lemma}\label{classisom} Every isometry of an ${\mathbb R}$-tree is either elliptic or
hyperbolic.\end{lemma} 
\demo Choose a point $x$ in $\Gamma$ and consider the intersection of $[x,\tau x]$ with
$\tau^{-1}\left([ x,\tau x]\right)$. This is a geodesic $[x,y]$. If $y=\tau y$ then $\tau$ is
elliptic. Otherwise, from the definition of $y$ it follows that the intersection of $[y,\tau y]$ with
$\tau^{-1}[y,\tau y]$ is the single point $y$, because any other point of intersection would also
give an extra point of intersection of $[x,\tau x]$ with $\tau^{-1}\left([ x,\tau x]\right)$. Thus
$\ell=\bigcup_n \tau^n[y,\tau y]$ is isometric to ${\mathbb R}$ and invariant under $\tau$. Now
$[x,\tau x]=[x,y]\cup[y,\tau y]\cup[\tau y,\tau x]$ and so $d(x,\tau x)\ge d(y,\tau y)$ with
equality if and only if $x$ is on $\ell$. Thus the points moved the minimal distance by $\tau$ are
exactly those points on $\ell$. Hence $\tau$ has axis $\ell$ and is thus hyperbolic\end{proof}

\begin{lemma}\label{translengthformula} Suppose that $\alpha,\beta$ are isometries of an ${\mathbb
R}$-tree and that ${\rm core}(\alpha)$ is disjoint from ${\rm core}(\beta)$. Then
$t(\alpha\beta)=t(\alpha)+t(\beta)+2d({\rm core}(\alpha),{\rm core}(\beta))$\end{lemma} \demo Let $a$ be a
point in ${\rm core}(\alpha)$ and $b$ a point in ${\rm core}(\beta)$ such that $d(a,b)$ is minimal. To see
that such points exist, choose points $a'$ in ${\rm core}(\alpha)$ and $b'$ in ${\rm core}(\beta)$. Define
$a$ by  $[a',a]=[a',b']\cap {\rm core}(\alpha)$ and $b$ by  $[b,b']=[a',b']\cap {\rm core}(\beta)$. Given any
arc, $\gamma,$ with one endpoint in ${\rm core}(\alpha)$ and the other endpoint in ${\rm core}(\beta)$ then
$\gamma$ contains $[a,b]$. For otherwise, since cores are convex we can join to $\gamma$ an arc in
${\rm core}(\alpha)$ connecting $a$ to $\gamma$ and another arc in ${\rm core}(\beta)$ connecting $b$ to
$\gamma$. Observe that the arcs joined to $\gamma$ meet $[a,b]$ only at $a$ or $b$. The union of
these three arcs is a path from $a$ to $b$ which contains a subset which is an embedded arc from $a$
to $b$. But $[a,b]$ is the unique arc in the ${\mathbb R}$-tree with these endpoints. Hence $\gamma$
contains $[a,b]$ as asserted. Hence ${\rm core}(\gamma)\ge d(a,b)$. Thus $a,b$ have the claimed property.

Define
$$\ell=[\alpha(b),\beta^{-1}(b)]=\alpha([b,a])\cup[\alpha(a),a]\cup[a,b]\cup[b,\beta^{-1}(b)]$$
Observe that $[\alpha(a),a]$ is contained in the ${\rm core}(\alpha)$ and $[b,\beta^{-1}(b)]$ is
contained in the ${\rm core}(\beta)$. Also $[a,b]$ intersects ${\rm core}(\alpha)\cup {\rm core}(\beta)$ only at $a,b$.
It follows from this that the four arcs on the right hand side only intersect at endpoints. Since
$a\in {\rm core}(\alpha)$ we have $t(\alpha)=d(a,\alpha(a))$ and similarly
$t(b)=d(\beta(b),b)=d(b,\beta^{-1}(b))$. Thus  $$\begin{array}{rcl}
{\rm core}(\ell)&=&d(\alpha(b),\alpha(a))+d(\alpha(a),a)+d(a,b)+d(b,\beta^{-1}(b))
\\ &=&d(\alpha(b),\alpha(a))+t(a)+d(a,b)+t(b)\\ 
&=&t(a)+2d(a,b)+t(b)\\
&=&t(a)+2d({\rm core}(\alpha),{\rm core}(\beta))+t(b).
\end{array}$$
We claim that $\ell\cap \alpha\beta(\ell)= \alpha(b),$ from
which it follows that  $\ell$ is contained in the axis of $\alpha\beta$ and is a
fundamental domain for $\alpha\beta$. The computation of ${\rm core}(\ell)$ above then completes the
proof.

Now $[b,\beta^{-1}(b)]$ is contained in ${\rm core}(\beta)$. Since ${\rm core}(\beta)$ is
invariant under $\beta$ it follows that $\beta([b,\beta^{-1}(b)])$ is contained in ${\rm core}(\beta)$.
Since $[b,a]\cap {\rm core}(\beta)=b$ it follows that $[b,a]\cap\beta([b,\beta^{-1}(b)])=b$. 
If the claim is false then 
$\ell\cap\alpha\beta(\ell)$ is an interval of positive length containing $\alpha(b)$. Now
$\alpha(b)$ is the start of $\ell$ and $\alpha([b,a])$ is a subinterval at the start of $\ell$.
Therefore this interval at the start of $\ell$ intersects the interval
$\alpha\beta([b,\beta^{-1}(b)])$ at the end of $\alpha\beta(\ell)$ in an interval of positive
length. Thus $\alpha([b,a])\cap\alpha\beta([b,\beta^{-1}(b)])$ has positive length. But this
implies $[b,a]\cap\beta([b,\beta^{-1}(b)])=b$ has positive length, a contradiction\end{proof}

\noindent An action of $G$ on an ${\mathbb R}$-tree $\Gamma$ is {\em minimal} if there is no  proper
subtree that is preserved by $G$. The action is {\em irreducible} if
there are two elements $a,b$ of $G$ which are both hyperbolic and whose translation axes are
disjoint. A {\em global fixed point} is a point in $\Gamma$ that is fixed by every element of $G$.

It can be shown that {\em irreducible} is equivalent to the statement that there is no point in
$\Gamma\cup\partial\Gamma$ which is fixed by all of $G$. Here, $\partial\Gamma$ is the {\em
boundary} of the ${\mathbb R}$-tree $\Gamma$. This boundary be described as either the set of ends of
$\Gamma,$ or as the Gromov-boundary of  $\Gamma$,
but we will not need this. If two hyperbolics $\alpha,\beta,$ have axes whose
intersection is compact, then for large $n$ the axes of $\alpha$ and $\beta^n\alpha\beta^{-n}$ are
disjoint.

\begin{lemma}\label{fixedpoint} 
Suppose $G$ is a finitely generated group of isometries of
an ${\mathbb R}$-tree $\Gamma$.
If every element  $G$ is elliptic, then there is a global fixed point.\end{lemma} \demo If $\alpha$ and $\beta$ are elliptic then
$$t(\alpha\beta)=2d({\rm core}(\alpha),{\rm core}(\beta))$$ Thus if $\alpha\beta$ is also elliptic then
${\rm core}(\alpha)$ intersects ${\rm core}(\beta)$. It follows that the intersection of the cores of a
finite number of elements which generate the group is non-empty and each point in the
intersection is fixed by every element in the group\end{proof}

It is clear that the action of $G$ on $K$ discussed above induces an action
$\s:G\longrightarrow \Isom(\Gamma)$ on $\Gamma$ by isometries.

\begin{addendum}\label{Gaction} Assume the hypotheses of (\ref{blowupgivesRtree}). Then $G$
acts isometrically on the ${\mathbb R}$-tree $\Gamma$ with no global fixed point. Moreover
 the action contains a hyperbolic. \end{addendum}
\demo Suppose that $G$ fixes the point $x$ in $\Gamma$. Then given $\e>0,$ for $n$ sufficiently
large, there is a point $x'_n$ in $(K,\tilde d_n)$ with image $x$ in $\Gamma$ such that $x'_n$ is
moved a $\tilde{d}_n$-distance at most $\e$ by all elements of $S$. Taking $\epsilon=1/2$  this
means that for all sufficiently large $n$ and all $\alpha\in S$ that:
$$d_\Hcal(\io_nx'_n , \alpha\io_nx'_n)\le\lambda_n/2$$ Thus $r_S(\iota_nx_n')\le\lambda_n/2$
which contradicts the definition of $\lambda_n$ and of the center of $\rho_n$. The conclusion
that some element of $G$ acts as a hyperbolic follows from lemma (\ref{fixedpoint}). \end{proof}

\begin{lemma}\label{invtsubtree} Suppose that $G$ acts by isometries on an ${\mathbb R}$-tree
$\Gamma$. Assume that some element of $G$ is hyperbolic. Let $\Gamma_-$ be the intersection
of all connected subsets of $\Gamma$ containing the axis of every hyperbolic. Then $\Gamma_-$
is the unique $G$-invariant minimal subtree of $\Gamma$. \end{lemma}
\demo Let $g$ be a hyperbolic element of $G$ with axis $\ell$. Then $t(g)>0$. Suppose that
$\Gamma'$ is a $G$-invariant subtree of $\Gamma$ then the restriction of $g$ to $\Gamma'$ is
an isometry with translation length in $\Gamma'$ at least as large as $t(g)$. By
(\ref{classisom}) it follows that $g|\Gamma'$ is a hyperbolic therefore it has an axis,
$\ell',$ which is invariant under $g|\Gamma'$. It is easy to check that the only line preserved by a
hyperbolic is its axis. Hence $\ell'=\ell$. Hence every invariant subtree contains the axis
of every hyperbolic. Therefore $\Gamma'$ contains $\Gamma_-$. Also $\Gamma_-$ is clearly
$G$-invariant and connected therefore an ${\mathbb R}$-tree\end{proof}

\begin{proposition}\label{integersimplicial} Suppose $G$ be a finitely generated group which acts by isometries on
an $\mathbb R$-tree $\Gamma$. Suppose the action and is minimal, and irreducible, and that the translation length of every element of $G$ is an integer. 
Then $\Gamma$ is a simplicial tree.\end{proposition} 
\demo We call a point $v$ in $\Gamma$ a {\em vertex} if $\Gamma-v$ has more than $2$ components.
We will show that the distance between vertices in $\Gamma$ is half of an integer. An inductive
argument now implies that $\Gamma$ is a simplicial tree.

Given two hyperbolic isometries $a,b$ in $G$ define $d$ to be the distance between their axes.
If $d>0$ then by (\ref{translengthformula}) we have $t(ab)=t(a)+t(b)+2d$. Since all
translation lengths are assumed to be integers this implies that $d$ is a half integer. 

Given two distinct vertices $u_1,u_2$ in $\Gamma$ we claim that, for $i=1,2,$ there are
hyperbolics $\tau_i$ with axis $\ell_i$ such that $\ell_i\cap[u_1,u_2]=u_i$. Then $[u_1,u_2]$
is the shortest arc between $\ell_1$ and $\ell_2$ and therefore by the previous remark
$d(u_1,u_2)$ is half an integer.

Let $P$ be the component of $\Gamma-interior[u_1,u_2]$ containing $u_1$. Since $u_1$ is a vertex
it follows that $P-u_1$ is not connected. Let $Q,R$ be the closures of two distinct components of
$P-u_1$. Since the action of $G$ is irreducible, there are hyperbolics $\alpha,\beta$ with
disjoint axes. Our subclaim is that there is a hyperbolic, $\tau_Q,$ with axis, $\ell_Q,$
contained in $Q-u_1$. Let $A$ be the axis of $\alpha$ and $B$ the axis of $\beta$. If the orbit of
$A$ under $G$ is disjoint from $Q$ then there is an invariant subtree given by taking the
intersection of all subtrees containing the orbit of $A$. This subtree does not contain $Q$
contradicting minimality. Thus we may assume that $A$ intersects $Q$. Then $A\cap Q$ is an
infinite half-line, since the axis $A$ can only enter or leave $Q$ at the point $u_1$. We may replace
$\alpha$ by $\alpha^{-1}$ if necessary to ensure that $\alpha(A\cap Q)\subset A\cap Q$. The axis of
$\alpha^n\beta\alpha^{-n}$ is $\alpha^n(B)$. For $n$ large $\alpha^n(B)$ contains points in $Q$
and does not contain $u_1$ therefore it is contained in $Q-u_1$. This proves the subclaim. Similarly
there is a hyperbolic, $\tau_R,$ with axis, $\ell_R,$ contained in $R-u_1$.

Suppose that $\delta$ is an arc connecting $\ell_Q$ to $\ell_R$. Then $interior(\delta)$ contains
$u_1$ because $Q\cap R=u_1$. Furthermore $\delta$ is contained in $Q\cup R$. Hence $\delta$
intersects $[u_1,u_2]$ only at $u_1$. Choose $\delta$ to be the shortest arc connecting the axes
of $\tau_Q$ and $\tau_R$. The proof of (\ref{translengthformula}) shows that the axis,
$\ell_1,$ of the hyperbolic $\tau_Q\tau_R$ contains $\delta$. Therefore $\ell_1$ intersects
$[u_1,u_2]$ in the single point $u_1$. Otherwise the intersection is an interval which contains
a subinterval of positive length in $\delta$ and this is impossible. Similarly one obtains a
hyperbolic with axis intersecting $[u_1,u_2]$ in the single point $u_2$. This proves the claim
\end{proof} 

\section{Length Functions and Compactifications.}

In this section we show how spaces of representations can be compactified by mapping the
representation space into the space of length functions on the group and compactifying this
latter space projectively.
 
Using the ideas in the previous section, any action, $\rho,$ of a group $G$ by
isometries on a geodesic space $\Hcal$ and any choice of a point, $x,$ in $\Hcal$ gives an
equivariant straight map $$\iota:K(G)\longrightarrow\Hcal$$ Let $d$ be the
pseudo-metric on $K(G)$ obtained by pull-back of the $d_\Hcal$ metric. This gives an
action, $\sigma,$ of $G$ on $K(G)$ by $d$-isometries. Given $g$ in $G$ the translation
length $t(\sigma g)$ is at least as large as the translation length $t(\rho g)$. It may be
larger, for example if $\rho(g)$ is hyperbolic and no point on the axis of $\rho(g)$ in $\Hcal$
is in  ${\rm image}(\iota )$. However there is a bound on the difference in the translation lengths
of $\rho g$ and $\sigma g$ which depends only on the thin triangles constant $\Delta$. We
now prove this.

\begin{lemma}\label{approxlength} Suppose that $\Hcal$ has $\Delta$-thin triangles and suppose
that $\tau$ is an isometry of $\Hcal$ and $x$ is any point in $\Hcal$. Then there is a
point $y$ in $[x,\tau x]$ such that $$|d_\Hcal(y,\tau y)-t(\tau)|\le 28\Delta$$\end{lemma}
\demo If $d(x,\tau x)=t(\tau)$ then $y=x$ works. Otherwise choose $z$ in $\Hcal$ such that 
$$d(z,\tau z)\ <\ \min\left(d(x,\tau x),t(\tau)+\Delta\right)$$
Every point in the orbit $\{ \tau^nz:n\in\ZZ \}$ of $z$ under the group
generated by $\tau$ is moved by $\tau$ the same distance that $z$ is moved. Thus we may assume
that $z$ is chosen in this orbit so that it is almost a closest point in the orbit to $x$ in the
sense that $d(z,x)<d(\tau^n z,x)+\Delta$ for all $n$. 

The first case is that there are points $y$ in $[x,\tau x]$ and $w$ in $[z,\tau
z]$ such that $d(y,w)\le2\Delta$. Now 
$$d(w,\tau w)\le d(w,\tau z)+ d(\tau z,\tau w) = d(w,\tau z)+d(z,w) = d(z,\tau z)$$ 
Also $$d(y,\tau y)\le d(y,w)+d(w,\tau w)+d(\tau w,\tau y) = 2d(y,w)+d(w,\tau w)$$
Combining these gives
$$d(y,\tau y)\le 4\Delta+d(z,\tau z) \le 4\Delta+t(\tau)+\Delta$$ This proves that this $y$
works in this case.

Next we prove that if the first case does not hold then $d(z,\tau z)\le20\Delta$. Assume
the contrary. Let $a,b$ be the points on $[z,\tau z]$ with $d(z,a)=d(b,\tau z)=5\Delta$. Let $a'$ be
the point on $[x,z]$ with $d(z,a')=5\Delta$ then we claim $d(a,a')\le4\Delta$. In the thin
quadrilateral with vertices  $x,\tau x, \tau z, z$ we have that every point on $[z,\tau z]$ is
within a distance of $2\Delta$ of the union of the two sides $[x,z]$ and $[\tau x,\tau z]$. This is
because the first case does not hold thus no point on $[z,\tau z]$ is within $2\Delta$ of $[x,\tau
x]$. In particular there is some point, $v$ say, in $[z,x]$ or
$[\tau z,\tau x]$ with $d(a,v)\le 2\Delta$. 
First suppose that
$v$ is in $[\tau z,\tau x]$ then:    $$\begin{array}{rcl} d(v,\tau z) &\ge& d(a,\tau z)-d(a,v)\\ &=&
d(z,\tau z)-d(z,a)-d(a,v)\\ &\ge& 20\Delta-5\Delta-2\Delta\\ &=& 13\Delta\end{array}$$ 
Hence $$\begin{array}{rcl}
d(\tau^{-1}z,x)&=&d(z,\tau x)\\
&\le& d(z,a)+d(a,\tau x)\\
&\le& 5\Delta+d(a,v)+d(v,\tau x)\\
&\le& 5\Delta+2\Delta+\left(d(\tau z,\tau x)-d(v,\tau z)\right)\\
&\le& 7\Delta+d(z,x)-13\Delta\\
&=&d(z,x)-6\Delta.
\end{array}$$ This implies that $d(z,x)>d(\tau^{-1}z,x)+\Delta$ which contradicts the choice of
$z,$ proving the subclaim that $v$ is in $[z,x]$.
 Since $d(a,z)=5\Delta$ and $d(a,v)\le2\Delta$ it follows that
$3\Delta\le d(v,z)\le 7\Delta$. Since $v$ and $a'$ are both on $[z,x]$ and $d(a',z)=5\Delta$ it
follows that $d(v,a')\le2\Delta$. Hence $d(a,a')\le4\Delta$ as claimed. Now define $b'=\tau a'$ thus
$b'$ is the point on $[\tau x,\tau z]$ with $d(\tau z,b')=5\Delta$. A similar argument shows that
$d(b,b')\le4\Delta$. It now follows that $$\begin{array}{rcl} d(a',\tau a')=d(a',b') &\le&
d(a',a)+d(a,b)+d(b,b')\\ & \le & 4\Delta+(d(z,\tau z)-d(z,a)-d(\tau z,b))+4\Delta\\ & \le & 8\Delta
+ (d(z,\tau z)-10\Delta)\\ & \le & d(z,\tau z)-2\Delta\\
& < & (t(\tau)+\Delta)-2\Delta\\
& < & t(\tau)\end{array}$$ It is impossible that $\tau$ moves $a'$ less
than $t(\tau)$ thus $d(z,\tau z)\le20\Delta$. It remains to show $y$ exists with this hypothesis.

Since the first case does not hold it follows that every
point on $[x,\tau x]$ is within $2\Delta$ of the union of $[x,z]$ and $[\tau x,\tau z]$. Hence
there is some point $y$ on $[x,\tau x]$ which is within a distance of $2\Delta$ of a point $c$
on $[x,z]$ and also of a point $e$ on $[\tau x,\tau z]$. This is proved by an easy continuity
argument. Now
$$\begin{array}{rcl}
|\ d(e,\tau z)-d(c,z)\ | & \le & d(e,c)+d(\tau z,z)\\
& \le & d(e,y)+d(y,c)+20\Delta\\
& \le & 24\Delta.
\end{array}$$ 
This choice of $y$ works because:
$$\begin{array}{rcl}
d(y,\tau y)& \le & d(y,e)+d(e,\tau c)+d(\tau c,\tau y)\\
& = & d(y,e)+d(e,\tau c)+d(c,y)\\
& \le & 2\Delta + d(e,\tau c) + 2\Delta\\
& = & 4\Delta + |d(e,\tau z)-d(\tau c,\tau z)|\\
& = & 4\Delta + |d(e,\tau z)-d(c,z)|\\
& \le & 4\Delta+24\Delta.
\end{array}$$ 
Observe that we use that $e,\tau z, \tau c$ all are on $[\tau x,\tau z]$ to obtain the second
equality.\end{proof}

\begin{corollary}\label{translength} Suppose that $G$ is a group and $K$ is the complete graph
with vertex set $G$. Suppose that $\Hcal$ is a  geodesic space with $\Delta$-thin triangles. The
action of $G$ on itself by left multiplication extends to an action, $\sigma,$ on $K$.  

Suppose that
$\rho:G\longrightarrow\Isom(\Hcal)$,
and $\iota:K \longrightarrow {\mathcal
H}$ is a $\rho$-equivariant straight map. Then
 $\sigma$ preserves $d=\iota^*d_{\Hcal}$ . Moreover, for
all $g$ in $G$ $$|t(\sigma g)-t(\rho g)|\le 28\Delta$$\end{corollary}
\demo Given $g$ in $G$ since $\io$ is an isometry onto its image it follows that $t(\rho g)\le
t(\sigma g)$. Consider the edge $[1,g]$ in $K$. Then $\iota([1,g])=[x,\rho(g).x]$ where
$x=\iota(1)$. The lemma (\ref{approxlength}) gives a point $y=\iota(z)$ such that
$$t(\sigma g) \le d_{K}(z,\sigma(g).z)=d_\Hcal(y,\rho(g).y)\le t(\rho g)+28\Delta$$\end{proof}

A {\it length function } on a group $G$ is a function $L:G\longrightarrow \mathbb R$. The set
of length functions on $G$ is denoted $\LF(G)$ and equals $\prod_{g\in G}{\mathbb R}$ with
the product topology. Thus a sequence $L_n$ in $\LF(G)$ converges to $L$ if for all
$g$ in $G$ we have that $L_n(g)$ converges to $L(g)$.

It is clear that conjugate isometries of $X$ have the same translation length, thus some
authors regard a length function as a real valued map from {\em  conjugacy classes.}
Given a representation $\rho:G\longrightarrow \Isom(X)$ we define the {\it length function }
$L_{\rho}:G\longrightarrow\mathbb R$ associated to $\rho$ by $$L_{\rho}(g)=t(\rho g)$$ for $g$ in
$G$. This in turn defines a map $$L:\Hom(G,\Isom(X))\longrightarrow\LF(G)$$
by $L(\rho)=L_{\rho}$.  We will consider the cases that $X$ is $\Hcal$ or an
${\mathbb R}$-tree.

In the previous section we introduced rescaling factors to obtain limiting ${\mathbb
R}$-trees. Thus we wish to discuss convergence of length functions in a projective sense. Note
that ${\mathbb R}-0$ acts on $\LF(G)$ by multiplication. The space of {\em projectivized length
functions} on a group $G,$ denoted $\PLF(G)$ is the set of projective equivalence classes of
non-zero length functions $$\PLF(G)\equiv\frac{\LF(G)-0}{{\mathbb R}-0}$$ This is given
the quotient topology of the subspace topology on $\LF(G)-0$. Thus a sequence of length
functions $L_n$ define projective classes which converge to the projective class of a length function
$L$ if there are $\lambda_n\ne0$ such that for all $g$ in $G$ $$\lim\ \lambda_n L_n(g)\ =\ L(g)$$ We
assume here that none of these length functions is identically zero. 

The following is the main result of this section:

\begin{proposition}\label{lengthlimit} Suppose that $G$ is a finitely generated group and
$\Hcal_n$ is a geodesic space with $\Delta$-thin triangles. Suppose that
$\rho_n:G\longrightarrow \Isom(\Hcal_n)$ is a sequence of representations which blows up. Let
$K$ be the complete graph with vertex set $G$. Let $x_n$ be an approximate center of $\rho_n$ and 
$\iota_n:K\longrightarrow\Hcal_n$ an equivariant straight map such that
$\iota_n(g)=\rho_n(g)x_n$, and set $d_n=\iota_n^*d_{\mathcal{H}_n}$.

Then there is a subsequence, still denoted $\{\rho_n\},$ and
$\lambda_n\to 0$ so that $\lambda_n\cdot d_n$ converges pointwise to a pseudo-metric $\od$ on $K$.
Furthermore, the metric space associated to $(K,\od)$ is an ${\mathbb R}$-tree $\Gamma$. 

The action
of $G$ on $K$ covers an action $\sigma:G\longrightarrow \Isom(\Gamma)$ and some element of $\sigma(G)$ is hyperbolic. Furthermore $\lim_{n\to\infty}
L(\rho_n)=L(\s)$ in $\PLF(G)$. 

Suppose, in addition, that $\sigma$ is irreducible and
that $L(\s):G\to\ZZ$. Then the minimal $G$-invariant subtree of $\Gamma$ is a simplicial
tree on which $G$ acts with the same length function as $\sigma$\end{proposition} \demo The
existence of the ${\mathbb R}$-tree $\Gamma$ comes from (\ref{blowupgivesRtree}) and (\ref{Gaction})
gives the action, $\sigma,$ of $G$ on $\Gamma$ and guarantees that for some element $g$ of $G$ that
$\sigma g$ is hyperbolic. Thus $t(\sigma g)>0$ so $L(\sigma)\ne0$ defines a (non-zero !) projective
length function. Now (\ref{translength}) implies that  $L(\rho_n)\to L(\s)$ in $\PLF(G)$.  To
see this, let $$\sigma_n:G\longrightarrow \Isom(K,\lambda_nd_n)$$ be the action action of $G$ on $K,$
with $K$ regarded as a pseudo-metric space using the rescaled pseudo-metric $\lambda_nd_n$. Then
$|t(\sigma_ng)-\lambda_nt(\rho_ng)|\le\lambda_n28\Delta$. Now $\lambda_n28\Delta\to0$ and so
$L(\sigma_n)-\lambda_nL(\rho_n)\to0$ in $\LF(G)$. Now $L(\sigma_n)\to L(\sigma)$ in
$\LF(G)$ hence in $\PLF(G)$. Also $\lambda_nL(\rho_n)$ has the same limit in
$\PLF(G)$ as $L(\rho_n)$. This proves the claim. If, in addition, $\sigma$ is irreducible and
$L(\sigma)$ is integer valued, then (\ref{integersimplicial}) supplies the $G$-invariant simplicial
sub-tree\end{proof}

We can embed $\LF(G)$ into
$$\overline{\LF(G)}\equiv\LF(G)\cup\PLF(G)$$ 
topologized so that a sequence $L_n\in\LF(G)$ converges to $[L]\in\PLF(G)$ if and only if
there are $\lambda_n\to0$ such that $\lambda_nL_n\to L$ in $\LF(G)$. 
It was shown in \cite{CS1} that if a sequence of representations of $G$ into $\SL_2({\CC})$ does
not blow up then there is a subsequence which, after suitable conjugacy, converges to a
representation of $G$ into $\SL_2{\CC}$. See \cite{CL} corollary (2.1) for a geometric proof of
this fact. 
 Then (\ref{lengthlimit})
gives:

\begin{corollary}\label{compactification} Suppose $G$ is a finitely generated group, and $\Delta>0$  and
 $\Hom_{\Delta}(G)$ is the set of all homomorphisms of $G$ into $\Isom(X)$ for all geodesic
spaces $X$ with $\Delta$-thin triangles. Then the image of $L:\Hom_{\Delta}(G)\longrightarrow
\LF(G)$ has compact closure in $\overline{\LF(G)}$. Every point
in $\overline{\LF(G)}\cap \PLF(G)$ is the projectivized length function
of an action of $G$ on an $\mathbb R$-tree
 which has no global fixed point\end{corollary} 
\demo
Since $G$ is countable, every sequence  in $\Hom_{\Delta}(G)$ has a subsequence
whose  length functions converge in $\overline{\LF(G)}$.
 It follows from (\ref{lengthlimit}) that limit points in $\PLF(G)$ are projectivized length functions of actions on ${\mathbb R}$-trees.  \end{proof}

 Specializing (\ref{lengthlimit})  to this situation, if the sequence
blows up, then there is a subsequence which converges, in $\PLF(G),$ to an action of $G$ on an
${\mathbb R}$-tree:

\begin{theorem}\label{blowupgivestree} Suppose that $G$ is a finitely generated group and
$\rho_n:G\longrightarrow \SL_2({\CC})$ is a sequence of representations. Then either there is a
subsequence which converges in $\Hom(G,\SL_2{\CC})$, or the sequence blows up and there
is an $\RR$-tree $\Gamma$ and
an action $\sigma:G\longrightarrow\Isom(\Gamma)$ with the following properties. 
\nb There is no point of $\Gamma$ fixed by $G;$ equivalently $L(\s)$ is not identically zero.  \nb
There is a subsequence $\rho_{n_i}$ such that $\lim L(\rho_{n_i})=L(\sigma)$ in ${\PLF}(G)$.\\ If,
in addition, $L(\sigma)$ is integer valued, and $\sigma$ is irreducible, then $\Gamma$ may be chosen to
be a simplicial tree.  \end{theorem}

\section{Length Functions and Valuations}
In this section we show that given a curve $C\subset \Hom(G,\SL_2{\CC})$ satisfying certain mild
hypotheses, and given an end $\epsilon$ of $C$, then the rescaled length functions of  representations
approaching $\epsilon$  converge to an integer valued function. This is the translation length
function of a certain irreducible action on a simplicial tree.

The translation length of an element $A\in\SL_2{\CC}$ with large trace is approximately
twice the logarithm of the absolute value of its trace. Now trace is a polynomial function on the
representation variety, and thus we are concerned with the logarithm of a polynomial when the
variables are large. The logarithm of a polynomial of one complex variable is approximately the
degree of the polynomial times the logarithm of the variable when the variable is large. 

For a polynomial $p$ of several independent variables the asymptotic behaviour of  $\log|p|$ may
be rather complicated. However the  restriction of $\log|p|$ to a {\em complex curve} behaves in
this respect like a polynomial function of one variable. The notion of the degree of a polynomial
of one variable is replaced by evaluating a certain {\em valuation} on the polynomial. This
valuation gives the rate of growth of the polynomial on the curve as one goes to infinity along
the curve towards a particular topological end of the curve. 

We will now give a geometric interpretation of the valuations we need.
Let $C$ be a complex curve in ${\CC}^n$ and let $\pi:{\CC}^n\longrightarrow{\CC}$ be a
linear projection such that $z=(\pi|C)$ is proper. It is not hard to see that such a projection
always exists. Then $C$ is a smooth 2\--dimensional manifold except at finitely many points (see
\cite{Mu}) and the derivative of $z$ is non\--zero except at finitely many points. Let $d$ be the
degree (as a map) of $z,$ which is the same as the number of pre-images of a regular value since
$z$ is complex differentiable.

\begin{proposition}\label{limits} Let $C$ be a complex affine algebraic curve in ${\CC}^n$. Then
there is an integer $d>0$ with the following property. Let $y:{\CC}^n\longrightarrow {\CC}$ be a
polynomial and let $\gamma:[0,\infty)\longrightarrow C$ be a proper arc. Then $\lim_{t\to\infty} \log|y(\gamma
t)|/\log|z(\gamma t)|=b/a$ for some integers $a,b$ with $|a|\le d$.\end{proposition}

\demo The closure of the image of $C$ in ${\CC}^2$ under the polynomial map $(y,z)$ is a
complex curve $K,$ (see  \cite{Mu}). A curve in ${\CC}^2$ is the zero set of a polynomial,
thus $K$ is the zero set of a polynomial $p(y,z)$. The {\it Newton polygon,} $\Newt(p),$ of
$p$ is the convex hull of the finite set of points in the plane: $$\{(m,n): \text{ the\ coefficient\ of
}\ y^mz^n\ \text{ in }\ p(y,z)\ \text{is\ not\ zero}\ \}$$ 
Since $z=(\pi|C)$ is proper, for $t$ large $|z(\gamma t)|$ is large. Thus, in looking at the order
of magnitude of a term in $p(y(\gamma t),z(\gamma t))$ with $t$ large, we may ignore the modulus
of the coefficients in $p$. 

The order of magnitude of $y^mz^n$ is $m\ \log|y|+n\
\log|z|\equiv\theta(m,n)$ which is a linear function of $(m,n)$. Thus for a given large $t$ the order of
magnitude of the term $a.y^mz^n$ in $p(y,z)$ is approximately given by the linear map $\theta(m,n)$. Now
$p(y(\gamma t),z(\gamma t))=0$ hence there cannot be a single term of $p(y,z)$ which is an order of magnitude
larger than all other terms in $p(y,z)$. The terms of approximately largest magnitude in $p(y,z)$
nearly lie along a level set for the linear function $\theta$. This implies that the terms of largest
order of magnitude are those terms lying along some edge, $e,$ of $\Newt(p),$ and that this edge is
almost contained in a level set of $\theta$.

 As $t\to\infty$ the level sets of $\theta$ become more
nearly parallel to the edge $e$. We claim that for all $t$ sufficiently large, the edge $e$ is the
independent of $t$. For otherwise, as $t$ increases and the edge $e$ changes to a new edge $e'$ and
there is a value of $t$ when the terms along both edge $e$ and edge $e'$ are the same order of
magnitude, and this is impossible.

Let $p_e(y,z)$ be the {\it edge polynomial } consisting of the sum of those terms of $p$ which
lie along $e$.  Now $p_e(y,z)=y^r z^s q(y^a z^b)$ where $a,b$ are coprime
integers such that $e$ has slope $b/a,$ and $q$ is a polynomial of one variable with $q(0)\ne0$.
Since the terms in $p_e$ are larger than any other terms in $p$ it follows that
$\lim_{t\to\infty}q(y^a(\gamma t)z^b(\gamma t))=0$. Hence $\lim_{t\to\infty}y^az^b=
\eta$ where $\eta$ is a (necessarily non-zero) root of $q$. Now $\log|y|,\log|z|\to\pm\infty$ along $\gamma$ so
that $$\lim_{t\to\infty} \log|y(\gamma t)|/\log|z(\gamma t)|=-b/a$$ For fixed large $z$ the equation $q(y^a
z^b)=0$ has $|a|$ distinct solutions, and thus $p(y,z)=0$ also has $|a|$ distinct solutions. But
since the degree of $z$ is $d$ we have $|a|\le d$ as claimed\end{proof}

\noindent Thus $\gamma$ determines a map $\ D_{\gamma}:\{rational\ functions\ on\ {\CC}^n\}
\longrightarrow\mathbb Q\ $ by $$D_{\gamma}(y)=\lim_{t\to\infty} \log|y(\gamma t)|/\log|z(\gamma t)|$$
Take the smallest integer $p$ such that $V_{\gamma}(y)=pD_{\gamma}(y)$ is an integer for all polynomials
$y$. Observe that such $p$ exists and that $p\le d!$ by the above. We claim that changing the arc
$\gamma$ by a proper homotopy does not change $V_{\gamma}$. Suppose that $\gamma'$ is properly homotopic to
$\gamma$. Then the argument in (\ref{limits}) shows that there is an edge $e$ for $\gamma$ containing the
terms of largest order of magnitude for large $t,$ and $D_{\gamma}(y)$ is determined by this edge.
Similarly there is an edge $e'$ for $\gamma'$. The proper homotopy gives a $1$-parameter family,
$\gamma_s,$ of proper arcs with $\gamma_0=\gamma$ and $\gamma_1=\gamma'$. Thus we get a $1$-parameter family of edges
$e_s$ for these arcs. Now all the edges in this family are the same, by the same reasoning in
(\ref{limits}) that the edge $e$ was indepenent of $t$. This proves the claim.

Thus the map $V_{\e}\equiv V_{\gamma}$ depends only on the end, $\e,$ of $C$ to which $\gamma$ converges.
Now ends of $C$ correspond to ideal points of the smooth model $\tilde{C}$ of the
projectivization of $C$ and all we have done is to construct the (negative of the) valuation
associated to an ideal point of $\tilde{C}$. Note that our sign convention means that
$V_{\e}(y)>0$ iff $\lim_{t\to\infty} |y(\gamma t)| =\infty$. We say that a function $f(t)$ {\em blows
up} as $t\to\infty$ if $\lim_{t\to\infty} |f(t)| =\infty$. Thus a rational function values
positive if it blows up as you go out to the end of the curve in question.

Let $G$ be a finitely generated group and $\Hom(G,\SL_2{\CC})$ the affine algebraic set of 
representations $\rho:G\longrightarrow \SL_2{\CC},$ following common usage we will call this the
{\it representation variety}. Let $C$ be a curve in $\Hom(G,\SL_2{\CC})$ and $\e$ an end of $C$ such
that there is some element $\alpha\in G$ for which ${\rm{trace}}(\rho\alpha)\to\infty$ as $\rho\to\e$. Let $V_{\e}$
be the associated valuation. It is an observation of Culler and Shalen that that if one works with a
curve of characters instead of representations, then every end has this property eg. by \cite{CL}
Corollary 2.1.

We next show that for a representation near an end, $\epsilon,$ of $C$ the translation lengths of
elements of $G$ are approximately certain integer multiples of a large parameter. The integer assigned
to an element $\alpha$ of $G$ is the twice the valuation associated to $\epsilon$ of the trace of
$\alpha$. The formula for translation length of an element $A$ of $\SL_2{\CC}$ is
$$t(A)=2\cosh^{-1}(|{\rm trace}(A)|/2)$$ For ${\rm{trace}}(A)$ large this is approximately $2\log|{\rm trace}(A)|$. Our
primary concern is with elements of $G$ whose translation length goes to infinity
as one goes out along $C$ towards $\epsilon$. Such elements act on the limiting tree with
non-zero translation length given by this integer. 

This is why, following Culler and Shalen, we introduce the function $$f_{\alpha}(\rho)=({\rm{trace}}\
\rho\alpha)^2-4\ $$ then $$\ t(\rho\alpha)/\log|f_{\alpha}(\rho)|\approx 1$$ when ${\rm{trace}}(\rho\alpha)$ is large.
If ${\rm{trace}}(A)\approx\pm 2$ this measures how fast the translation length goes to zero.

Let $\rho_n$ be a sequence converging to an end $\e$ of $C$. Using the above approximation to
$t(\rho_n\alpha)$ and proposition (\ref{limits}) there are $\lambda_n$ such that for all $\alpha$ in $G$
$$\lim_{n\to\infty}\lambda_nt(\rho_n\alpha)=\max\{V_{\e}(f_{\alpha}),0\}$$ The corresponding formula of
Morgan and Shalen \cite{MS} Proposition [II.3.15] has a negative sign because our $V$ is the
negative of the valuation that they use. By (\ref{blowupgivestree}) there is a subsequence which
converges to an action $\sigma:G\longrightarrow Aut(\Gamma)$ on an ${\mathbb R}$-tree $\Gamma$. 
Hence

\begin{proposition}\label{valuationgiveslengths} Let $G$ be a finitely generated group and
let $C$ be a an affine algebraic curve of representations of $G$ into $\SL_2{\CC}$.
Suppose that $\epsilon$ is an end of $C$ and that $\rho_n$ is a sequence in $C$ which blows
up and converges to $\epsilon$. Define a length function on $G$ by
$p(g)=\max\{V_{\e}(f_{g}),0\}$ for each $g$ in $G$. Then there is an action, $\sigma,$ of $G$ on an
${\mathbb R}$-tree such that 
$\lim_{\rho\to\epsilon}L(\rho)=L(\sigma)=L(p)$ in $\PLF(G)$\end{proposition}
\demo By (\ref{lengthlimit}) there is a subsequence $\rho_{n_i}$ and an action $\sigma$ of $G$ on
an ${\mathbb R}$-tree such that $\lim_{i\to\infty}L(\rho_{n_i})=L(\sigma)$ in $\PLF(G)$. From the
above discussion $\lim_{i\to\infty}L(\rho_{n_i})=\lim_{\rho\to\epsilon}L(\rho)=L(p)$\end{proof}

In particular this means that we may scale metric on the $\mathbb R$\--tree $(\Gamma,\od)$
obtained as a limit of the representations $\rho_n$ so that the translation length function
for it is integer valued. 

\begin{lemma}\label{irreducible} Suppose $G$ is a finitely generated group and $C$ is a connected component of a curve in $\Hom(G,\SL_2\CC)$ and
$\Gamma$ is an $\RR$-tree and suppose
\begin{itemize}
\item  $\rho_n$ is a sequence in $C$.
\item  $\alpha\in G$ and  ${\rm{trace}}(\rho_n\alpha)\to\infty$.
\item  $\exists\ \rho'\in C$ and $\rho'(G)$ is not virtually solvable.
\item  $\exists\ \sigma:G\longrightarrow\Isom(\Gamma)$ and $\sigma(\Gamma)$ has no global fixed point. 
\item $L(\rho_n)$ converge to $L(\sigma)$ in $\PLF(G)$.
\end{itemize} Then $\sigma$ is irreducible. \end{lemma}

\demo We may suppose that a subsequence of $\rho_n$ has been chosen which converges to some end
$\epsilon$ of $C$. Since $C$ contains a representation which is not virtually solvable  then for
$\rho$ near $\epsilon$ we have that $\rho$ is not virtually solvable. The Tits alternative says that
a subgroup of a linear group is either virtually solvable or else contains a non-abelian free
subgroup. (The proof of this for $\SL_2{\CC}$ is elementary.) Thus there is a conjugate, $\beta,$
of $\alpha$ such that the fixed points of $\rho\beta$ in $S^2_{\infty}$ are both distinct from those
of $\rho\alpha$. Conjugate $\rho$ so that $$A = \rho(\alpha) = \left( \begin{array}{cc} a & 0  \\ 0 &
1/a \end{array} \right)\; \qquad B= \rho(\beta) = \left( \begin{array}{cc} b+c & b^2-c^2-1  \\ 1 &
b-c \end{array} \right). $$ This may be done so that $a\to\infty$ as $\rho\to\epsilon$. Then 
$$t_p\equiv{\rm{trace}}(A^pBA^{-p}B^{-1}) = (1-b^2+c^2)(a^{2p}+a^{-2p}) + 2(b^2-c^2)$$ Observe that
${\rm{trace}}(\rho\beta)=2b$ thus $b$ may be regarded as a polynomial function on $C$.   Similarly
${\rm{trace}}(\rho\alpha\beta)=a(b+c)-(b-c)/a$ so that $c$ may also be regarded as a rational function on
$C$.

 We claim that for all sufficiently large $p$ that $t_p\to\infty$ near $\epsilon$. To
see this, note that if $(1-b^2+c^2)=0$ then $B$ is lower triangular hence fixes $0\in
S^2_{\infty}$. But $A$ fixes $0$ also, and this contradicts the choice of $\beta$.  Let
$V_{\epsilon}$ be the valuation assoicated to $\epsilon$. Given two functions, one of which grows
faster than the other, then the rate of growth of their sum equals the rate of growth of the
larger one, thus $$V_{\epsilon}(f)>V_{\epsilon}(g) \quad\implies\quad V_{\epsilon}(f+g)=
V_{\epsilon}(f)$$ For $p$ sufficiently large we have 
$$V_{\epsilon}((1-b^2+c^2)(a^{2p}+a^{-2p}))=V_{\epsilon}(1-b^2+c^2) +
V_{\epsilon}(a^{2p}+a^{-2p}) > V_{\epsilon}(2(b^2-c^2))$$ This is because $a\to\infty$ and so 
$V_{\epsilon}(a^{2p}+a^{-2p})\ge 2p$. Thus for $p$ large, $V_{\epsilon}(t_p)>0$ so $t_p\to\infty$
near $\epsilon$ as claimed.

It follows that $\sigma\alpha^p$ and $\sigma\beta$ are hyperbolics and that their commutator is also
hyperbolic. Thus the axes of these two elements intersect in at most a compact interval. Then
$\sigma\alpha^p$ and a conjugate of it by any sufficiently large power of $\sigma\beta$ have
disjoint axes. Hence $\sigma$ is irreducible\end{proof}

In summary we have the main theorem about compactifying curves of representations:

\begin{theorem}\label{theorem}\label{mainthm} Let $C$ be a  connected component of an algebraic  curve in the representation variety $\Hom(G,\SL_2{\CC})$ such that
\nb  $\exists\ \rho\in C$ such that $\rho(G)$ is not solvable.
\nb $\exists\ \rho,\rho'\in C$ that are not conjugate.\\
If $\e$ is an end of $C$ then there is a metric simplicial tree $T$, and an irreducible action  $\sigma:G\longrightarrow\Isom(T)$ with the following property. Suppose that $\rho_n$
is any sequence in $C$ converging to $\epsilon$. Then $\lim L(\rho_n)=L(\sigma)$ in
$\PLF(G)$\end{theorem}
We remark that if all the representations on $C$ are virtually solvable then it is
easy to show that there is an action $\sigma$ of $G$ by isometries on ${\mathbb R}$ such that
$\lim L(\rho_n)=L(\sigma)$.

\bibliographystyle{amsplain}

\end{document}